\documentclass[reqno,12pt]{amsart}

\usepackage{amssymb,amsmath,amsthm,amsfonts,euscript}
\usepackage{graphicx}
\usepackage{pgf}
\usepackage{tikz}
\usetikzlibrary{cd}

\usepackage[colorinlistoftodos]{todonotes}

\usepackage[top=1in,bottom=1in,left=1in,right=1in]{geometry}

\newcommand{\C}{\mathbb{C}}
\newcommand{\K}{\mathbb{K}}
\newcommand{\Q}{\mathbb{Q}}

\newcommand{\Z}{\mathbb{Z}}

\newcommand{\la}{\lambda}

\newcommand{\tH}{\tilde{H}}

\newcommand{\Hilb}{\mathrm{Hilb}}
\newcommand{\loc}{\mathrm{loc}}

\newcommand{\Utor}{\mathbf{U}^{\mathrm{tor}}}

\newcommand{\aff}{\mathrm{aff}}

\newcommand{\gl}{\mathfrak{gl}}
\renewcommand{\sl}{\mathfrak{sl}}
\newcommand{\hsl}{\widehat{\mathfrak{sl}}}

\newcommand{\Y}{\mathbb{Y}}

\newcommand{\core}{\mathrm{core}}
\newcommand{\quot}{\mathrm{quot}}

\newcommand{\kb}{\overline{\kappa}}

\newcommand{\Mac}{\mathsf{M}}
\newcommand{\sT}{\mathsf{T}}
\newcommand{\Pol}{\mathbf{Pol}}
\newcommand{\Sym}{\Lambda}
\newcommand{\fS}{\mathfrak{S}}

\newcommand{\sx}{\mathsf{X}}

\newcommand{\sk}{\mathsf{k}}

\newcommand{\Ju}{\underline{J}}

\newcommand{\Iu}{\underline{I}}
\newcommand{\sA}{\mathsf{A}}

\newcommand{\id}{\mathrm{id}}
\newcommand{\tquot}{\widetilde{\quot}}
\newcommand{\Xd}{X^\bullet}
\newcommand{\mud}{{\mu^\bullet}}
\newcommand{\cP}{\mathcal{P}}
\newcommand{\La}{\Lambda}

\newcommand{\trianglerightneq}{\mathrel{\ooalign{\raisebox{-0.5ex}{\reflectbox{\rotatebox{90}{$\nshortmid$}}}\cr$\triangleright$\cr}\mkern-3mu}}
\newcommand{\triangleleftneq}{\mathrel{\reflectbox{$\trianglerightneq$}}}
\newcommand{\pair}[1]{\langle \, #1\, \rangle}

\newtheorem{lem}{Lemma}
\newtheorem{thm}[lem]{Theorem}

\newtheorem{dfn}[lem]{Definition}

\theoremstyle{remark}
\newtheorem{rem}[lem]{Remark}
\newtheorem*{ex}{Example}

\numberwithin{equation}{section}
\numberwithin{lem}{section}

\begin{document}

%[Difference Operators for wreath Macdonald polynomials and applications]

\title[Difference Operators for wreath Macdonald polynomials]{Difference Operators for wreath Macdonald polynomials}
\author{Daniel Orr}
\address{}
\email{dorr@vt.edu}
\author{Mark Shimozono}
\email{mshimo@math.vt.edu}
\address{
	Department of Mathematics (MC 0123),
	460 McBryde Hall, Virginia Tech,
	225 Stanger St.,
	Blacksburg, VA 24061 USA}

\maketitle

\begin{abstract}
We give explicit $q$-difference operators acting diagonally on wreath Macdonald $P$-polynomials in finitely many variables.
\end{abstract}

\section{Introduction}

\textit{Wreath Macdonald polynomials} \cite{Hai} originate from the geometry of $\Gamma$-fixed loci in Hilbert schemes of points in the plane, where $\Gamma$ is a finite cyclic subgroup of $SL(2,\C)$. Unlike the ordinary Macdonald polynomials (which are recovered when $\Gamma$ is trivial), their very \textit{existence} relies on the use of extremely deep geometric methods \cite{BezF,BezK}.

A recent result of Wen \cite{Wen} provides a somewhat more direct characterization of wreath Macdonald polynomials. This result is based on the \textit{quantum toroidal algebra} \cite{GKV} of type $\mathfrak{gl}_r$, where $r=|\Gamma|$, together with a remarkable symmetry of this algebra known as Miki's automorphism \cite{Mi}. However, Wen's result is still quite far from an explicit characterization of wreath Macdonald polynomials. 
One encounters serious technical obstacles when trying to use it to extract concrete information about wreath Macdonald polynomials, of which very little is in fact known.

The aim of this note is to introduce a new family of $q$-difference operators for wreath Macdonald polynomials, which are a {\em direct} generalization of Macdonald's original $q$-difference operators for ordinary Macdonald polynomials. In particular, we provide explicit first-order $q$-difference operators $\{\Mac^{(i)}\}_{i\in I}$ (where $I=\Z/r\Z$) which act diagonally on wreath Macdonald $P$-polynomials in finitely many variables with an $I$-graded version of the ordinary Macdonald operator eigenvalues (see Theorem~\ref{T:main}). Our operators are sufficient to uniquely determine these polynomials and their lifts to multi-symmetric functions.

We have a complete proof that our operators act diagonally on wreath Macdonald polynomials, but we postpone it to \cite{OSW}, where we will also construct higher order wreath Macdonald difference operators. Our approach is based on: (i) Wen's result mentioned above, and (ii) Negut's shuffle realization of the quantum toroidal algebra \cite{Neg}. Based on the rich applicability of Macdonald's original difference operators, we anticipate a number of applications of the operators $\{\Mac^{(i)}\}_{i\in I}$, which we will take up in future papers.

\section*{Acknowledgements}

We thank Mark Haiman and Joshua Wen for helpful discussions. D. O. was partially supported by a Collaboration Grant for Mathematicians from the Simons Foundation.

\section{Macdonald polynomials}\label{S:mac}

Let $\fS_N$ be the symmetric group and consider the space $\Sym_N=\K[x_1,\dotsc,x_N]^{\fS_N}$ of symmetric polynomials in variables $x_1,\dotsc,x_N$ with coefficients in $\K=\Q(q,t)$. Let $\Y$ be the set of all integer partitions $\lambda=(\lambda_1,\lambda_2,\dotsc)$ where $\lambda_1\ge \lambda_2\ge\dotsm\ge 0$, and let $\Y_N$ be the set of all integer partitions $\la=(\la_1,\dotsc,\la_N)$ with at most $N$ parts.

\subsection{Origins}
The {\em Macdonald operator} is the $q$-difference operator
\begin{align}\label{E:M}
\Mac=\sum_{k=1}^N\prod_{\substack{l=1\\l\neq k}}^N\frac{tx_k-x_l}{x_k-x_l}\cdot\sT_{q,x_k}
\end{align}
where the rational functions appearing in $\Mac$ are understood as multiplication operators and
\begin{align}\label{E:T}
\sT_{q,x_k}f(x_1,\dotsc,x_N) &= f(x_1,\dotsc,qx_k,\dotsc,x_N).
\end{align}
The operator $\Mac$ stabilizes $\Sym_N$, i.e., the denominators clear and the result is symmetric. The importance of this operator stems from the fact that it belongs to a quantum integrable system:
\begin{thm}[\cite{Mac}]\label{T:M}
\begin{enumerate}
\item
There exist explicit, algebraically independent $q$-difference operators $\Mac=\Mac_1, \Mac_2, \dotsc, \Mac_N,$ 
satisfying $[\Mac_k,\Mac_l]=0$ for all $1\le k,l\le N$.

\item The operators $\Mac_1,\dotsc,\Mac_N$ have a joint eigenbasis $\{P_\la\}$ indexed by $\la\in\Y_N$. The symmetric polynomials $P_\la$ are uniquely characterized as their joint eigenfunctions, together with the normalization
\begin{align*}
P_\lambda \in m_\lambda + \bigoplus_{\mu\triangleleftneq\la} \K m_\mu
\end{align*}
where $m_\lambda$ denotes a monomial symmetric polynomial and $\unlhd$ is the dominance order on partitions (restricted to $\Y_N$).
\end{enumerate}
\end{thm}

%\begin{align}
%\Mac_d=t^{d(d-1)/2}\sum_{\substack{I\subseteq\{1,\dotsc,n\}\\|I|=d}}\prod_{\substack{i\in I\\j\notin I}}\frac{tx_i-x_j}{x_i-x_j}\prod_{i\in I}\sT_{q,x_i}
%\end{align}

The symmetric polynomials $P_\la$ are the {\em Macdonald polynomials} (of type $GL_N$). We will call $\Mac_2,\dotsc,\Mac_N$ the {\em higher Macdonald operators}. The eigenvalues of $\Mac$ are:
\begin{align}\label{E:M-e}
\Mac\cdot P_\lambda &= e(q,t) P_\lambda,\qquad 
e(q,t)=\sum_{k=1}^N q^{\lambda_k} t^{N-k}.
\end{align}
The eigenvalues of $\Mac_k$ are given by substituting the summands of $e(q,t)$ into the elementary symmetric function $e_k$.

%Maybe work in orthogonality.

%One has
%\begin{align}
%\Mac_a P_\la &= e_a(q^{\la_1}t^{N-1},q^{\la_2}t^{N-2},\dotsc,q^{\la_N})P_\la
%\end{align}
%for all $1\le a\le N$, where $e_a$ denotes the $a^{\mathrm{th}}$ elementary symmetric function.

The Macdonald polynomials $P_\lambda$ can be stabilized into elements of the $\K$-algebra $\Lambda$ of symmetric functions in infinitely many variables \cite{Mac}. One also considers the modified Macdonald symmetric functions $\tH_\la$, which are obtained from $P_\lambda$ by an explicit transformation (see, e.g., \cite{Hai}).

\subsection{Geometry of Hilbert schemes}

%A series of conjectures, made by Macdonald, were a major driving force in the study of Macdonald polynomials. One of the most famous of these, the {\em positivity conjecture}, was proved by Mark Haiman.

The modified Macdonald symmetric functions $\tH_\la$ are significant due to their positivity when expanded in Schur symmetric functions, which was the content of Macdonald's famous \textit{positivity conjecture}. Haiman's proof \cite{Hai} of this conjecture recast it as a statement about the geometry of Hilbert schemes $\Hilb_m=\Hilb_m(\C^2)$ of points in the plane. The deepest part of Haiman's proof concerned the existence and properties of a rank $m!$ vector bundle $P_m$ on $\Hilb_m$ known as the Procesi bundle.

It is convenient to work with the union $\Hilb=\bigsqcup_{m\ge 0}\Hilb_m$ of all Hilbert schemes. Let us consider their localized $T$-equivariant $K$-groups:
\begin{align*}
K_T(\Hilb)_\loc=\bigoplus_{m\ge 0}K_T(\Hilb_m)_\loc
\end{align*}
%of equivariant $K$-groups (of coherent sheaves), 
where $T=(\C^\times)^2$ acts naturally on $\C^2$ and hence on $\Hilb$, and $\loc$ stands for localization with respect to equivariant scalars $R(T)\cong\Z[q^{\pm 1},t^{\pm 1}]$, which is achieved by tensoring with $\K\cong\mathrm{Frac}(R(T))$. The points of $\Hilb$ can be identified with ideals $I\subset \C[x,y]$ of finite codimension. By the localization theorem, $K_T(\Hilb)_\loc$ has a $\K$-basis given by the classes of $\{[I_\lambda]\}$ of $T$-fixed points in $\Hilb$, which are precisely the monomial ideals $I_\lambda=(y^{\lambda_1},xy^{\lambda_2},x^2y^{\lambda_3},\dotsm)$ where $\lambda\in\Y$ is any partition.

\begin{thm}[\cite{Hai:pos}]\label{T:H}
The Procesi bundles $P_m$ induce an isomorphism of $\K$-vector spaces
\begin{align}\label{E:H-iso}
K_T(\Hilb)_\loc\cong\Lambda
%K_{\fS_m\times T}(\C^{2m})_\loc=
\end{align}
under which a fixed point class $[I_\la]$ is sent to the modified Macdonald function symmetric function $\tH_\la$. This isomorphism realizes $\tH_\la$ as the bigraded $\fS_m$-character of the fiber $P_m|_{I_\la}$, where $|\la|=m$, implying the Macdonald positivity conjecture.
\end{thm}

\section{Wreath Macdonald polynomials}

%Now we come to the wreath setting. 
For a fixed integer $r>1$, we consider the finite cyclic subgroup 
$$\Gamma = \Bigg\{\begin{pmatrix}\zeta & \\ & \zeta^{-1}\end{pmatrix} : \zeta^r=1\Bigg\} 
\subset SL_2(\C). $$
Haiman \cite{Hai} proposed a generalization of the {\em geometric} theory of Macdonald polynomials, replacing $\Hilb$ by the $\Gamma$-fixed locus $\Hilb^\Gamma$ consisting of finite-codimensional ideals $I\subset\C[x,y]$ which are $\Gamma$-invariant.\footnote{In fact, Haiman's generalization applies conjecturally to arbitrary finite subgroup $\Gamma\subset SL_2(\C)$; however, we only consider cyclic $\Gamma$, as very little is known outside of this case.} %Before we can describe this more precisely, we need to establish some notation.
Let $\chi : \Gamma \to \C^\times$ be the representation of $\Gamma$ given by
\begin{align}
\chi:\begin{pmatrix}\zeta & \\ & \zeta^{-1}\end{pmatrix}\mapsto\zeta
\end{align}
The irreducible representations of $\Gamma$ are given by powers of $\chi$: $\mathrm{Irr}(\Gamma)=\{\chi^i : i\in I\}$ where $I=\Z/r\Z$. The McKay correspondence identifies $I$ with the vertices of the Dynkin diagram of type $A^{(1)}_{r-1}$. Let $Q_\aff=\oplus_{i\in I}\Z\alpha_i$ be the root lattice of the affine Lie algebra $\hsl_r$, $Q_\aff^+=\oplus_{i\in I}\Z_{\ge 0}\alpha_i$ the positive affine root cone, and $\delta=\sum_{i\in I}\alpha_i$ the minimal imaginary root. %, and we identify $Q=Q_\aff/\Z\delta$ with the root lattice of $\sl_r$.
Then the fixed locus $\Hilb^\Gamma$ is the disjoint union of connected components $\Hilb^\Gamma_\alpha$ indexed by $\alpha=\sum_{i\in I} m_i\alpha_i\in Q_\aff^+$, namely\footnote{The name ``wreath'' originates from the fact that $\Hilb^\Gamma_\alpha$ is a crepant resolution of $\C^{2n}/\Gamma_n$ for some $n$ (determined by $\alpha$), where $\Gamma_n=\Gamma^n\rtimes \fS_n$ is a wreath product group.}:
\begin{align*}
\Hilb^\Gamma_\alpha = \{I\in \Hilb^\Gamma : \text{$R/I$ has $\Gamma$-character $\sum_{i\in I} m_i\chi^i$}\}.
\end{align*}
$\Hilb^\Gamma$ and $\Hilb$ share the same $T$-fixed points $\{I_\lambda\}_{\lambda\in\Y}$, since $\Gamma\subset T$.

\subsection{Haiman's conjecture}
In \cite{Hai} Haiman conjectured the following generalization of Theorem~\ref{T:H}: there exists a wreath Procesi bundle on each component of $\Hilb^\Gamma$ and that together these afford an isomorphism
\begin{align}\label{E:wreath-iso}
K_T(\Hilb^\Gamma)_\loc \cong \Lambda^\Gamma := \bigoplus_{\gamma\in Q} \Lambda^{\otimes I}
\end{align}
with properties generalizing Macdonald positivity, where $Q=Q_\aff/\Z\delta$ is the root lattice of $\sl_r$. Parallel to Theorem~\ref{T:H}, the images of fixed point classes $\{[I_\lambda]\}$ under this isomorphism are the modified {\em wreath Macdonald symmetric functions} $\tH^\Gamma_\la$, whose precise definition will be recalled below.

\begin{rem}
In the isomorphism \eqref{E:wreath-iso}, the element $\gamma\in Q$ records the $r$-core of a partition $\lambda\in\Y$ indexing a fixed-point class on the left-hand side. The factors in the tensor power $\Lambda^{\otimes I}$ of symmetric functions have to do with the $r$-quotient of $\lambda\in\Y$, under the core-quotient decomposition \cite{JK}.
\end{rem}

Haiman's conjecture was proved by Bezrukavnikov and Finkelberg \cite{BezF}. Their proof relies on deep, non-explicit results of Bezrukavnikov and Kaledin \cite{BezK} involving quantization in positive characteristic.

\subsection{Definition of wreath Macdonald polynomials}
\label{SS:wreath macs} The following definitions of $\tH^\Gamma_\la$ and their monic forms $P^\Gamma_\la$ are due to Haiman \cite{Hai}, with existence guaranteed by \cite{BezF}.

\subsubsection{$\tH$-version}
For $\la\in\Y$ denote by $\core_r(\la)\in\Y$ and $\quot_r(\la)\in\Y^I$ the classical $r$-core and $r$-quotient \cite{JK}. We identify $\la\in\Y$ with its diagram $D(\la)=\{(a,b)\in\Z_{\ge0}^2\mid 0 \le a < \la_{b+1} \}.$
We identify the set of all $r$-core partitions with the root lattice $Q$ of $\sl_r$ as follows. 
Let $\kappa:\Y\to Q_\aff$ be the map $\kappa(\la) =  \sum_{(a,b)\in\la} \alpha_{b-a}$ with subscripts taken modulo $r$ as usual, and let $\kb = \mathrm{cl} \circ \kappa: \Y\to Q$ be the composition of $\kappa$ with the projection $\mathrm{cl}:Q_\aff\to Q$. It is known that $\kb$ restricts to a bijection from the set of $r$-core partitions to $Q$ \cite{Hai}. 

Let $\trianglerighteq^\Gamma$ denote the following restriction of the standard dominance order $\trianglerighteq$ on $\Y$ \cite{Mac}: two partitions may only be compared if they have the same $r$-core and the same number of boxes. %\footnote{According to \cite{Hai:unpublished} one may consider Nakajima quiver varieties involving another parameter $u\in \fS_r$ which uses a more general stability condition. We will require this later.} 
The map $\quot_r$ gives a bijection $\Y^\gamma\cong  \Y^I$ from the set 
$\Y^\gamma$ of partitions with fixed $r$-core $\gamma\in Q$, to the set $\Y^I$ of $r$-tuples of partitions
\cite[\S 6.2]{G}.

For $\la\in\Y$ we write $s_\la[\Xd]\in \Lambda^{\otimes I}$ for the tensor Schur function $s_\mud[\Xd] = \prod_{i\in I} s_{\mu^{(i)}}[X^{(i)}]$ indexed by $\mud=\tquot_r(\la)\in\Y^I$ where $\tquot_r(\la)$ is the reversal of the $r$-tuple $\quot_r(\la)$. (Here $X^{(i)}$ stands for an alphabet of variables $x^{(i)}_1, x^{(i)}_2, \dotsc$ corresponding to the $i$-th factor in $\Lambda^{\otimes I}$.)

For $a\in \K$ we denote by $\cP_{\id- a \chi^{-1}}$ the $\K$-algebra automorphism of $\La^{\otimes  I}$ which maps the power sum $p_d[X^{(i)}]$ of degree $d$ living at vertex $i\in I$ of the McKay graph of $\Gamma$ as follows: $$\cP_{\id- a \chi^{-1}}(p_d[X^{(i)}])= p_d[X^{(i)}] - a^d\, p_d[X^{(i-1)}].$$ Let $\pair{\cdot,\cdot}$ be the tensor product of Hall pairings on $\La^{\otimes I}$, with respect to which the tensor Schur functions are orthonormal.

The wreath Macdonald polynomials $\tH^\Gamma_\la\in\La^{\otimes I}$ are uniquely defined by \cite{Hai,BezF}
\begin{enumerate}
    \item $\cP_{\id - q \chi^{-1}}(\tH^\Gamma_\la) \in \K^\times s_\la[\Xd] + \bigoplus_{\mu\trianglerightneq^\Gamma \la} \K s_\mu[\Xd]$.
    \item $\cP_{\id - t^{-1} \chi^{-1}}(\tH^\Gamma_\la) \in \K^\times s_\la[\Xd] + \bigoplus_{\mu\triangleleftneq^\Gamma \la} \K s_\mu[\Xd]$.
    \item $\pair{\tH^\Gamma_\la, s_n[X^{(0)}]} = 1$ where $n$ is the total number of boxes in $\tquot_r(\la)$.
\end{enumerate}
\begin{rem} \label{R:basis}
For each fixed $r$-core $\gamma$, the $\tH^\Gamma_\la$ for $\la\in\Y^\gamma$ form a $\K$-basis of $\La^{\otimes I}$. So do the $s_\la[\Xd]$ for $\la\in \Y^\gamma$; for each $\gamma$ this is the same tensor Schur basis of $\La^{\otimes I}$ but indexed (and ordered) differently.
\end{rem}

\subsubsection{$P$-version}
Define $J_\la^\Gamma\in \La^{\otimes I}$ by 
\begin{align}\label{E:J def}
J_\la^\Gamma = \mathcal{P}_{\id - t^{-1}\chi^{-1}} (\tH_\la^\Gamma).
\end{align}
Then $P_\la^\Gamma\in\La^{\otimes I}$ is by definition the unique $\K^\times$-multiple of $J_\la^\Gamma$ for which the coefficient of $s_\la[\Xd]$ is $1$. Due to the triangularity implied by the definition of $\tH^\Gamma_\la$, the $P_\la^\Gamma$ form a $\K$-basis of $\La^{\otimes I}$ as $\la$ runs over $\Y^\gamma$ for any fixed $r$-core $\gamma$.
\begin{rem}\label{R:1}
When $r=1$, one recovers the ordinary Macdonald $P$ functions as $P^\Gamma_\la(q,t) = P_\la(q^{-1},t)$.
\end{rem}

\section{Wreath Macdonald difference operators}

\subsection{Symmetric functions and symmetric polynomials}
Our difference operators will characterize the images of $\{P^\Gamma_\lambda\}_{\lambda\in\Y}$
under projections to finitely many variables.
Let $N_\bullet=(N_i)_{i\in I}\in(\Z_{\ge 0})^I$ be a dimension vector and consider the polynomial ring
\begin{align*}
\Pol_{N_\bullet} &= \K[\{x^{(i)}_k : i \in I, 1\le k\le N_i\}]
\end{align*}
and its subring
\begin{align}
\Sym_{N_\bullet} = (\Pol_{N_\bullet})^{\fS_{N_\bullet}}
\end{align}
consisting of polynomials which are invariant with respect to the product of symmetric groups 
$\fS_{N_\bullet}=\prod_{i\in I} \fS_{N_i}$,
with the factor $\fS_{N_i}$  permuting the variables $\{x^{(i)}_1,\dotsc,x^{(i)}_{N_i}\}$, for each $i\in I$. We write $N=\sum_{i\in I}N_i$ for the total number of variables. 

We write
$
\pi_{N_\bullet} : \Lambda^{\otimes I} \to \Sym_{N_\bullet}
$
for the natural projection, which is given on simple tensors $f=\otimes_{i\in I} f_i$ %(where each $f_i$ belongs to $\Lambda$) 
by
\begin{align}
\pi_{N_\bullet}(f)=\prod_{i\in I} f_i(x^{(i)}_1,\dotsc,x^{(i)}_{N_i}).
\end{align}
For $f\in\Lambda^{\otimes I}$ we use the shorthand $f(X_{N_\bullet})$ to denote $\pi_{N_\bullet}(f)$.

For $\gamma\in Q$ indexing a summand of $\Lambda^\Gamma$, a dimension vector $N_\bullet\in(\Z_{\ge 0})^I$ will be called \textit{$\gamma$-compatible} if it satisfies
\begin{align}\label{E:compat}
N_{i} - N_{i-1} &= \langle \alpha_i^\vee,-\gamma\rangle, \qquad \text{for each $i\in I$.}
%= \langle \alpha_i,-w_0\kb(\lambda)\rangle
\end{align}
where the subscripts in $I=\Z/r\Z$ are understood in a cyclic sense and $\alpha^\vee_0$ means the classical projection of the $0$-th affine simple coroot. 

In what follows, we have to assume that $N_\bullet$ is $\gamma$-compatible. But there is no loss of generality in this assumption, since \eqref{E:compat} continues to hold if we add a constant to each $N_i$. Thus we can make the components of $N_\bullet$ as large as needed to fully capture the $P_\lambda^\Gamma$ as symmetric functions.

\subsection{Difference operators}
Now we are ready to introduce our wreath Macdonald difference operators. % $\{\Mac^{(i)}\}_{i\in I}$. 
The Macdonald operator $\Mac$ in \eqref{E:M} is a sum over all choices of a variable $x_k \ (1\le k\le N)$. Our operators $\Mac^{(i)}$ are expressed as sums over all nonempty subsets of the variables $\bigcup_{i\in I}\{x^{(i)}_1,\dotsc,x^{(i)}_{N_i}\}$ with at most one variable selected at each vertex. We record these subsets as pairs $\Ju=(J,\sk)$, where $J$ is a nonempty subset of $I$ and $\sk$ is an integer-valued function on $J$ such that $1\le \sk(j)\le N_j$ for each $j\in J$; the variables chosen to be in our subset are $\{x^{(j)}_{\sk(j)}\}_{j\in J}$. 

Let $\mathcal{S}$ be the set of all such pairs and consider some $\Ju=(J,\sk)\in \mathcal{S}$. For any $i\in I$, we define:
\begin{align*}%\label{E:sx-all}
\text{$\sx_{\Ju}^{(i)}=q^{m}x_{\sk(i+m)}^{(i+m)}$,\quad where $m\ge 0$ is minimal such that $i+m\in J$.}
\end{align*}
For $j\in J$ we have 
$\sx_{\Ju}^{(j)}=x^{(j)}_{\sk(j)}$, which is simply the variable at vertex $j$ selected by $\Ju$. The remaining $\sx_{\Ju}^{(i)}$ for $i\in I\setminus J$ are obtained by cyclically propagating the selected variables and multiplying by a factor of $q$ at each step.

\begin{ex}
Let $I=\Z/3\Z$ and $N_\bullet=(2,2,2)$. An element of $\Pol_{N_\bullet}$ is a polynomial
$$f=f(x^{(0)}_1,x^{(0)}_2\mid x^{(1)}_1,x^{(1)}_2 \mid x^{(2)}_1,x^{(2)}_2);$$ dividers separate the  variables at each vertex. Let $\Ju=(J,\sk)$ where $J=\{0,1\}$ and $\sk(0)=\sk(1)=1$.
That is, the selected subset is $\{x^{(0)}_1,x^{(1)}_1\}$. Then
\begin{align*}
\sx_{\Ju}^{(0)} &= x^{(0)}_1,\quad
\sx_{\Ju}^{(1)} = x^{(1)}_1,\quad
\sx_{\Ju}^{(2)} = qx^{(0)}_1,
\end{align*}
where $x^{(0)}_1$ is cyclically propagated to the left and mutlplied by $q$ to determine $\sx_{\Ju}^{(2)}$.
\end{ex}

\begin{dfn}
For any $i\in I$, define the wreath Macdonald difference operator
\begin{align}\label{E:Mi}
\Mac^{(i)}
&=\sum_{\Ju\in\mathcal{S}^{(i)}}(-1)^{|J|}\cdot \sA_{\Ju}^{(i)}\cdot \sT_{q,\Ju}
\end{align}
where $\mathcal{S}^{(i)}\subset \mathcal{S}$ is the subset determined by the condition $i-1\in J$,
\begin{align}
\mathsf{A}_{\Ju}^{(i)}
=
\frac{\sx_{\Ju}^{(i)}}{\sx_{\Ju}^{(0)}} &\cdot \prod_{j\in I}\frac{\sx_{\Ju}^{(j-1)}}{\sx_{\Ju}^{(j-1)}-t\sx_{\Ju}^{(j)}}
\cdot \prod_{j\in J\setminus\{i-1\}}\frac{\sx_{\Ju}^{(j+1)}}{\sx_{\Ju}^{(j+1)}-q^{-1}\sx_{\Ju}^{(j)}}\notag\\
\label{E:A-coeff}
&\cdot \prod_{j\in I\setminus J}\frac{\displaystyle\prod_{l=1}^{N_{j-1}}(t\sx_{\Ju}^{(j)}-x^{(j-1)}_l)}{\displaystyle\prod_{l=1}^{N_j}(\sx_{\Ju}^{(j)}-x^{(j)}_l)}
\cdot \prod_{j\in J}
\frac{\displaystyle\prod_{l=1}^{N_{j-1}}(t\sx_{\Ju}^{(j)}-x^{(j-1)}_l)}{\displaystyle \sx_{\Ju}^{(j)}\cdot\prod_{\substack{l=1\\l\neq \sk(j)}}^{N_j}(\sx_{\Ju}^{(j)}-x^{(j)}_l)},
\end{align}
and $\sT_{q,\Ju}$ is the operator on $\Pol_{N_\bullet}$ given by:
%\begin{align}\label{E:T-act}
%\sT_{\Ju}^* =\ &\text{\em for each $j\in J$, replace $\sx_{\Ju}^{(j)}=x^{(j)}_{\sk(j)}$ with $q^{-m}\sx_{\Ju}^{(j+m)}=q^{-m}x_{\sk(j+m)}^{(j+m)}$},\\
%&\text{\em  where $m>0$ is minimal such that $j+m\in J$.}\notag
%\end{align}
%We note that $\sT_{\Ju}^*$ can be equivalently described as follows:
\begin{align}\label{E:T-act-alt}
\sT_{q,\Ju} =\ &\text{\em for each $j\in J$, replace $\sx_{\Ju}^{(j)}$ with $q\sx_{\Ju}^{(j+1)}$.}
\end{align}
\end{dfn}

\begin{thm}\label{T:main}
For any $\gamma\in Q$, $\lambda\in\Y^\gamma$, and $\gamma$-compatible $N_\bullet\in(\Z_{\ge 0})^I$, one has\footnote{ $P_\lambda^\Gamma(X_{N_\bullet};q^{-1},t)$ means that we project to finitely many variables and send $q$ to $q^{-1}$; cf. Remark~\ref{R:1}.}
\begin{align}\label{E:eigenequation}
\Mac^{(i)}\cdot P_\lambda^\Gamma(X_{N_\bullet};q^{-1},t) &= e_\lambda^{(i)}(q,t)P_\lambda^\Gamma(X_{N_\bullet}; q^{-1},t).
\end{align}
where the eigenvalues $e_\lambda^{(i)}(q,t)\in\K$ are given by (recall that $\chi\in\mathrm{Irr}(\Gamma)$ satisfies $\chi^r=1$):
\begin{align}\label{E:eigenvalues}
\sum_{i\in I} e_\lambda^{(i)}(q,t)\chi^i = \sum_{k=1}^N q^{\lambda_k} t^{N-k} \chi^{k-\lambda_k}\in R(T\times\Gamma)\cong \Z[q^{\pm 1},t^{\pm 1},\chi^{\pm 1}].
\end{align}
\end{thm}

\begin{ex}
Continuing the example above, we have
%We have $i=1$, $\nd=(2,2,2)$, and $\Ju=(J,\sk)$ where $J=\{0,1\}$ and $\sk(0)=\sk(1)=1$.
\begin{align*}%\label{E:A-coeff-ex}
\mathsf{A}_{\Ju}^{(1)}=\frac{\sx_{\Ju}^{(1)}}{\sx_{\Ju}^{(0)}}
&\cdot \frac{\sx_{\Ju}^{(2)}}{\sx_{\Ju}^{(2)}-t\sx_{\Ju}^{(0)}}\frac{\sx_{\Ju}^{(0)}}{\sx_{\Ju}^{(0)}-t\sx_{\Ju}^{(1)}}\frac{\sx_{\Ju}^{(1)}}{\sx_{\Ju}^{(1)}-t\sx_{\Ju}^{(2)}}
\cdot \frac{\sx_{\Ju}^{(2)}}{\sx_{\Ju}^{(2)}-q^{-1}\sx_{\Ju}^{(1)}}\\
&\cdot \frac{t\sx_{\Ju}^{(2)}-x^{(1)}_1}{\sx_{\Ju}^{(2)}-x^{(2)}_1}\frac{t\sx_{\Ju}^{(2)}-x^{(1)}_2}{\sx_{\Ju}^{(2)}-x^{(2)}_2}
\cdot \frac{t\sx_{\Ju}^{(0)}-x^{(2)}_1}{\sx_{\Ju}^{(0)}}
\frac{t\sx_{\Ju}^{(0)}-x^{(2)}_2}{\sx_{\Ju}^{(0)}-x^{(0)}_2}
\cdot \frac{t\sx_{\Ju}^{(1)}-x^{(0)}_{1}}{\sx_{\Ju}^{(1)}}
\frac{t\sx_{\Ju}^{(1)}-x^{(0)}_2}{\sx_{\Ju}^{(1)}-x^{(1)}_2}\notag
\end{align*}
and
$\mathsf{T}_{q,\Ju}(f) = f(qx^{(1)}_1,x^{(0)}_2 \mid q^2x^{(0)}_1,x^{(1)}_2 \mid x^{(2)}_1,x^{(2)}_2)
$.
\end{ex}

We hope this example convinces the reader that not even the most clever guessing would lead from the ordinary Macdonald operator $\Mac$ to its wreath generalization $\Mac^{(i)}$. One really needs the quantum toroidal algebra $\Utor_r$ to find these operators. 

\subsection{Recovering ordinary Macdonald operators}

Let us take $\Gamma=1$. Then $I=\{0\}$ and $N_\bullet=(N_0)=(N)$. The only possibly subsets of variables are singletons $\{x^{(0)}_k\}$ where $1\le k\le N$. In terms of $\Ju=(J,\sk)$, we have $J=I$ and $\sk(0)=k$. The operator $\sT_{q,\Ju}$ sends $x^{(0)}_k$ to $qx^{(0)}_k$, which agrees with $\sT_{q,x^{(0)}_k}$ in \eqref{E:M}. Finally, we have
\begin{align*}
\mathsf{A}_{\Ju}^{(0)}
&=
\frac{1}{1-t}
\cdot 
\frac{\displaystyle\prod_{l=1}^{N}(tx^{(0)}_k-x^{(0)}_l)}{\displaystyle x^{(0)}_k\cdot\prod_{\substack{l=1\\l\neq k}}^{N}(x^{(0)}_k-x^{(0)}_l)}
=-\prod_{\substack{l=1\\l\neq k}}^N\frac{tx^{(0)}_k-x^{(0)}_l}{x^{(0)}_k-x^{(0)}_l}
\end{align*}
and we see that $\Mac^{(0)}$ reduces to $\Mac$ from \eqref{E:M}.

\subsection{Relation to Shoji's operators}
We also mention an interesting similarity between our operators and some introduced by Shoji \cite{Sho}.
Let us consider only the terms in \eqref{E:Mi} which have $J=I$, meaning that a variable is selected at each vertex in $I$. Then, for $\Iu=(I,\sk)$, we have $\sx_{\Iu}^{(j)}=x^{(j)}_{\sk(j)}$ for all $j\in I$ and
\begin{align}\label{E:A-coeff-I}
\mathsf{A}_{\Iu}^{(i)}
=
(-1)^r\frac{\sx_{\Iu}^{(i)}}{\sx_{\Iu}^{(0)}}
\prod_{j\in I\setminus\{i-1\}}\frac{\sx_{\Iu}^{(j+1)}}{\sx_{\Iu}^{(j+1)}-q^{-1}\sx_{\Iu}^{(j)}}
\cdot \prod_{j\in I}\prod_{\substack{l=1\\l\neq\sk(j)}}^{N_j}
\frac{t\sx_{\Iu}^{(j+1)}-x^{(j)}_l}{\sx_{\Iu}^{(j)}-x^{(j)}_l}.\notag
\end{align}
Shoji \cite{Sho} considered the operator
\begin{align}
\mathsf{S}=\sum_{\Iu=(I,\sk)}
\prod_{j\in I}\prod_{\substack{l=1\\l\neq\sk(j)}}^{N_j}
\frac{t\sx_{\Iu}^{(j+1)}-x^{(j)}_l}{\sx_{\Iu}^{(j)}-x^{(j)}_l}\cdot \sT_{q,\Iu}
\end{align}
and used it to characterize a different family of Macdonald-like polynomials associated with the wreath product groups $\Gamma_n$. One might expect these to relate to wreath Macdonald polynomials. Haiman observed some ``accidental'' agreement between the two \cite{Hai}. Our operators clarify the genuine difference between Shoji's polynomials and wreath Macdonald polynomials. The operators $\Mac^{(i)}$ for wreath Macdonald polynomials are substantially more complicated than Shoji's operator $\mathsf{S}$.

%Clarify the relationship to Shoji's operators and his Hall-Littlewood functions.

\subsection{The proof}

Let us give an indication of how we discovered Theorem~\ref{T:main} and how we will prove it and extend it to higher wreath Macdonald operators in \cite{OSW}. The operators $\{\Mac^{(i)}\}_{i\in I}$ arise from the {\em horizontal Heisenberg subalgebra}  in the quantum toroidal $\gl_r$ algebra $\Utor_r$ \cite{GKV} acting on the space $\Lambda^\Gamma$, which affords the vertex representation of \cite{Sai}.

In its natural presentation, $\Utor_r$ comes equipped with an obvious Heisenberg subalgebra, the {\em vertical Heisenberg subalgebra}. In both the vertex representation $\Lambda^\Gamma$ and the geometric representation $K_T(\Hilb^\Gamma)_\loc$ constructed in \cite{VV}, the vertical Heisenberg algebra acts in a straightforward way. In contrast, it is difficult to describe the horizontal Heisenberg subalgebra and its actions on $\Lambda^\Gamma$ and $K_T(\Hilb)_\loc$. Miki's automorphism \cite{Mi} swaps the vertical and horizontal Heisenberg subalgebras of $\Utor_r$. 

Extending a result of \cite{FT1,SV2} to the wreath setting, Wen \cite{Wen} proved that the isomorphism \eqref{E:wreath-iso} respects the actions of $\Utor_r$ on both sides, up to some nontrivial scalars and a twist by the Miki automorphism. 
Thus, Wen's result implies that the horizontal Heisenberg algebra acts in $\Lambda^\Gamma$ like the vertical Heisenberg algebra acts in $K_T(\Hilb^\Gamma)_\loc$. From this one immediately deduces that the horizontal Heisenberg algebra acts diagonally on wreath Macdonald polynomials. With some care, this leads to the right-hand side of the eigenoperator equation \eqref{E:eigenequation} in Theorem~\ref{T:main}.
 
Obtaining the left-hand side, i.e., the explicit formula for the wreath Macdonald operators $\Mac^{(i)}$, requires substantial additional work. Using some known Miki automorphism images from \cite{T}, we transport degree-one elements in the horizontal Heisenberg algebra to the shuffle algebra \cite{Neg} and use this to compute the action on symmetric polynomials by a very careful iterated residue computation and further cancellation. In the passage from symmetric functions in infinitely many variables to symmetric polynomials in finitely many variables, skewing operators must be converted to the $q$-translations $\sT_{q,\Ju}$ in \eqref{E:T-act-alt}. When the dust settles, one obtains the formula for $\Mac^{(i)}$ in \eqref{E:Mi}. 

Finally, we note that we have extensively verified Theorem~\ref{T:main} using software for wreath Macdonald polynomials developed in Sage by the second named author.

%The sophistication of the algebraic objects involved here explains why it would be impossible to guess Theorem~\ref{T:main} as a generalization of Theorem~\ref{T:M}.

\end{document}